\documentclass[11pt,a4paper]{article}
\usepackage[utf8]{inputenc}
\usepackage[american]{babel}      % language
\usepackage[T1]{fontenc}
\usepackage{amsmath,amssymb,amsthm}
\usepackage{enumitem}
\usepackage{graphicx}
\usepackage{xcolor}   % pacchetto per scrivere a colori
\usepackage[left=3cm,right=3cm,top=3cm,bottom=3cm]{geometry}
\usepackage{setspace}
\title{{\bf Modular convergence of Steklov sampling operators\\ in Orlicz spaces}}
         
\author{ {\bf Danilo Costarelli} \quad and  \quad {\bf Erika Russo} \\  
Department of Mathematics and Computer Science \\
            University of Perugia\\
        1, Via Vanvitelli, 06123 Perugia, Italy    \\  
{\small {\tt danilo.costarelli@unipg.it} - {erika.russo1@studenti.unipg.it}} }

\date{}

\newcommand{\miu}{\leq}

\newcommand{\N}{\mathbb{N}}
\newcommand{\R}{\mathbb{R}}

\newcommand{\be}{\begin{equation}}
\newcommand{\ee}{\end{equation}}

\newtheorem{definition}{Definition}[section]
\newtheorem{theorem}[definition]{Theorem}
\newtheorem{lemma}[definition]{Lemma}
\newtheorem{proposition}[definition]{Proposition}

\newtheorem{corollary}[definition]{Corollary}
\theoremstyle{remark}
\newtheorem{remark}[definition]{Remark}
\begin{document}

\maketitle 

\begin{abstract}
In this paper, we deal with the family of Steklov sampling operators in the general setting of Orlicz spaces. The main result of the paper is a modular convergence theorem established following a density approach. To do this, a Luxemburg norm convergence for the Steklov sampling series based on continuous functions with compact support, and a modular-type inequality in the case of functions in Orlicz spaces has been preliminary proved. As a particular case of general theory, the results in $L^p$, in the Zygmund (interpolation), and in the exponential spaces are deduced. A crucial aspect in the above results is the choice of both band- and duration- limited kernel functions satisfying the partition of the unit property; to provide such examples an equivalent condition based on the Poisson summation formula and the computation of the Fourier transform of the kernel has been employed. 
\vskip0.3cm
\noindent
  {\footnotesize AMS 2010 Mathematics Subject Classification: 42A16; 47G30; 46E30; 41A05}
\vskip0.1cm
\noindent
  {\footnotesize Key words and phrases: Steklov averages, Steklov sampling operators, modular convergence, Orlicz spaces, Fourier transform, Poisson summation formula} 
\end{abstract}

\section{Introduction} \label{sec1}

The theory of Steklov sampling operators has been introduced and studied in \cite{3} in the $L^p$ setting, $1 \miu p < +\infty$, and in the spaces of bounded, continuous and uniformly continuous signals in $\mathbb{R}$. The main idea behind the introduction of this new family of sampling-type operators resides in the fact that, thanks to the use of high-order Steklov-type integrals, one can expect to improve the approximation performances of the already known families of sampling-type series. This is the main advantage that has been widely discussed (and obviously proved) in \cite{3}. 

As mentioned above, the idea consists in reconstructing a given signal $f$ using a sampling series based on a kernel function $\chi \colon \mathbb{R} \to \mathbb{R}$, which is a discrete approximate identity, and using net of sample values in a certain specific form $f_{r,w}(k/w)$; namely they are Steklov integrals of order $r$ evaluated at the nodes $\frac{k}{w}$, $k \in \mathbb{Z}$, $w>0$, according with the definition proposed by Sendov and Popov in \cite{1} and introduced to give a constructive proof of the Brudnii density theorem (\cite{25}).

Steklov means are commonly used in Approximation Theory and Fourier Analysis to provide regular approximations of not necessarily regular functions. In the context of signal processing, this method is particularly suitable because real-world signals often suffer from uncertainties such as time-jiitter and roundoff errors.

In this work, we extend the results of \cite{3} to the more general framework of Orlicz spaces; the main result here established is a modular convergence theorem, which has been proved following a density-type approach. To do this, a Luxemburg norm convergence for the Steklov sampling series based on continuous functions with compact support, and a modular-type inequality in the case of functions in Orlicz spaces has been preliminary proved. In particular, from the latter inequality it is also possible to deduce that the Steklov sampling series belong to a given Orlicz space $L^\varphi(\R)$ provided that $f$ itself belongs to $L^\varphi(\R)$. 

One of the main advantages of having approximation results at the disposal in this very general context is that Orlicz spaces represent a generalization of several known and widely used functional spaces. For instance, for suitable choices of the functions $\varphi$ (that must be a so-called $\varphi$-function), we are able to find the already mentioned Lebesgue spaces, the Zygmund (or interpolation) spaces, which are very useful in the theory of partial differential equations, and the exponential spaces that are employed for embedding theorems between Sobolev spaces.

A crucial aspect of the above results is represented by the choice of the kernel function $\chi$ defining the Steklov sampling operators. As a function $\chi$ we can consider both band- and duration- limited kernels satisfying the partition of the unit property, together with some other non-restrictive assumptions. In particular, concerning the partition of the unit property, in order to check if it is satisfied or not, one can resort to the application of a suitable equivalent condition based on the Poisson summation formula and the computation of the Fourier transform of the involved function. Examples of kernels satisfying such conditions are the Fejér kernel, the Jackson type kernel, the central B-splines, and several others (see, e.g., \cite{CORO,ACAR,DRA,COR,OZER}).

Finally, for the sake of completeness, we recall that the origin of the classical theory of sampling is due to the celebrated Whittaker-Kotel'nikov-Shannon (WKS) sampling theorem, which provides a fundamental interpolation formula for signals that are simultaneously band-limited and with finite energy. The reconstruction formula is provided by a cardinal series based on a family of uniform spaced sample values and the classical $sinc$-function (see, e.g., \cite{15, 16, 17}). The above theorem, even if it is very elegant from a mathematical point of view, presents some disadvantages from the point of view of applications, due to some consequences of the required assumptions on the reconstructed signal.

These considerations motivated the introduction of the theory of sampling operators. Among the most studied sampling-type operators we can mention the generalized sampling (GS) operators, introduced by the German mathematician P. L. Butzer (see \cite{24, 16, 17}) . GS operators are defined by means of a sequence of sample values of the form $f\left(\frac{k}{w}\right)$, $k \in \mathbb{Z}$, $w>0$, which become very suitable for reconstructing signals in continuous function spaces.

These operators in the context $L^p$ present some issues; for these reasons, in \cite{2} the Kantorovich sampling operators (KS) have been introduced, replacing in the definition of the GS operators the sample values $f\left(\frac{k}{w}\right)$ with certain averages of $f$ in suitable intervals containing $\frac{k}{w}$ (see also \cite{20,BAJ,fuzzy,CPV4}).

The KS operators were revealed to be more suitable than the GS operators to reconstruct not necessarily continuous signals, such as $L^p$ functions. For more details see \cite{2}.

The Steklov sampling operators represent the natural extension of the KS operators, in which the first-order averages of $f$ have been replaced by high-order Steklov integrals. 

The paper is structured as follows.

In Section \ref{sect2}, we first briefly discuss the theory of Orlicz spaces, introducing key concepts such as the Luxemburg norm, norm and modular convergence, together with some examples.

In Section \ref{sect3}, we introduce the theory of Steklov sampling operators and discuss the main properties of the considered kernel functions. 

In Section \ref{sect4}, we present the main approximation results established in the paper and described above.

While in Section \ref{Particular cases}, and Section \ref{sec6}, special spaces of Orlicz spaces and examples of kernels have been discussed. Finally, in Section \ref{sec7}, some final remarks and future research directions have been discussed.

%%%

\section{Orlicz spaces and modular convergence} \label{sect2}

Before delving into the development of this article, we introduce some notations that will be used in what follows.

We denote by $C(\mathbb{R})$ the space of all uniformly continuous and bounded functions on $\mathbb{R}$, endowed with the usual sup-norm $\| \boldsymbol{\cdot} \|_\infty$, and by $C_C(\mathbb{R}) \subset C(\mathbb{R})$ the subspace of the elements having compact support.

Moreover, $M(\mathbb{R})$ is the linear space of all (Lebesgue) measurable real
(or complex) functions defined on $\mathbb{R}$.

Now, we recall some basic notions about Orlicz spaces (see, e.g., \cite{MAL,HA2015}).

A real function $\varphi \colon \mathbb{R}_{0}^+\to\mathbb{R}_{0}^+$ is called a \emph{$\varphi$-function} if it is non-decreasing, continuous, and satisfies the following conditions: 
\begin{enumerate}
    \item $\varphi(0)=0$,
    \item  $\varphi(u)>0$, for $u>0$,
    \item $\varphi(u) \to +\infty$ as $u \to +\infty$.
\end{enumerate} 

For a fixed $\varphi$-function $\varphi$, we can consider the functional $I^\varphi \colon M(\mathbb{R})\to [0, +\infty]$, defined by 
\begin{equation}
 I^\varphi[f]:= \int_{\mathbb{R}} \varphi\left(| f(x)|\right) \ dx, \quad f \in M(\mathbb{R}).
\end{equation}
It is well-known that $I^\varphi$ is a modular on $M(\mathbb{R})$ generated by the $\varphi$-function $\varphi$, i.e., $I^\varphi$ satisfies the following assumptions:
\begin{enumerate}
    \item $I^\varphi[f]=0 \Leftrightarrow f=0$;
    \item $I^\varphi[f]=I^\varphi[-f]$;
    \item $I^\varphi[\alpha f+\beta g] \leq I^\varphi[f]+I^\varphi[g]$, for every $\alpha, \beta>0$ such that $\alpha+\beta=1$.
\end{enumerate}

The \emph{Orlicz space} generated by $\varphi$ is 
\begin{equation}
 L^\varphi (\mathbb{R}):= \left\{f \in M(\mathbb{R}) \colon \ I^\varphi [\lambda f]< + \infty \text{ for some } \lambda > 0 \right\}.
\end{equation}

It is clear that the Orlicz spaces are very general vector spaces that generalize, for instance, $L^p$-spaces. Indeed, it is not difficult to see that, setting $\varphi(u)=u^p$, $1 \leq p < +\infty$, it turns out that $L^\varphi(\mathbb{R})=L^p(\mathbb{R})$.

We can also introduce a notion of norm on $L^\varphi(\mathbb{R})$, called the \emph{Luxemburg norm}, that can be defined by \[ \| f \|_\varphi := \inf \left\{\lambda > 0: \ I^\varphi \left[\frac{f}{\lambda}\right] \leq \lambda \right\};\]
in case of $\varphi$ convex, the above notion reduces to
\be \label{Lux-norm}
 \| f \|_\varphi := \inf \left\{\lambda > 0: \ I^\varphi \left[\frac{f}{\lambda}\right] \leq 1 \right\}. 
\ee

From (\ref{Lux-norm}) it is possible to deduce a (strong) notion of convergence in $L^\varphi(\mathbb{R})$; a net of functions $(f_w)_{w>0} \subset L^\varphi(\mathbb{R})$ is \emph{norm convergent} to a function $f \in L^\varphi(\mathbb{R})$, i.e., $\|f_w -f\|_\varphi \to 0$ for $w \to +\infty$ if and only if $I^\varphi[\lambda(f_w-f)] \to 0$ as $w \to +\infty$, for every $\lambda>0$.

We can also introduce in $L^\varphi(\mathbb{R})$ a weak notion of convergence with respect to that induced by the Luxemburg norm: the modular convergence.
We will say that the net of functions $(f_w)_{w \in \mathbb{N}} \subset L^\varphi(\mathbb{R})$ is \emph{modularly convergent} to a function $f \in L^\varphi(\mathbb{R})$ if 
$$\lim_{k \to +\infty} I^\varphi[\lambda(f_k - f)]=0$$ for some $\lambda>0$.

Obviously, the Luxemburg norm convergence implies the modular convergence, while the converse implication is true if and only if the $\varphi$-function $\varphi$ satisfies the $\Delta_2$-condition (see, e.g., \cite{9, 11}). For the sake of completeness, we recall that $\varphi$ satisfies the $\Delta_2$-condition if there exists a constant $M>0$ such that
$$\varphi(2u) \leq M\varphi(u), \quad u \in \mathbb{R}_0^+.$$
The space
\begin{equation}
 E^\varphi(\mathbb{R}):= \left\{f \in M(\mathbb{R}): \ I^\varphi [\lambda f]< +\infty \text{ for every } \lambda > 0 \right\} 
\end{equation}
is a vector subspace of $L^\varphi(\mathbb{R})$ and it is called \emph{the space of all finite elements} of $L^\varphi (\mathbb{R})$. It is easy to show that $C_C(\mathbb{R}) \subset E^\varphi(\mathbb{R})$.

 In general, $E^\varphi (\mathbb{R})\subseteq L^\varphi(\mathbb{R})$, i.e., $E^\varphi (\mathbb{R})$ is a proper subspace of $L^\varphi(\mathbb{R})$. However, these two spaces coincide if and only if $\varphi$ satisfies the $\Delta_2$-condition.

 As a last property of Orlicz spaces, we recall that $C_C(\mathbb{R})$ is modularly dense in $L^\varphi(\mathbb{R})$ (see \cite{14}). 

We now give some remarkable examples of Orlicz spaces. For instance, we can mention two cases of $\varphi$-functions satisfying the $\Delta_2$-condition, which are the previously introduced $\varphi(u)=u^p$, for $1 \leq p < + \infty$ and $\varphi(u)=u^\alpha \log^\beta(e+u)$ for $\alpha \geq 1$ and $\beta>0$. 
In particular, the Orlicz spaces generated by the latter $\varphi$-functions are the spaces $L^\alpha \log^\beta L$, which are called the Zygmund (or interpolation) spaces.

In all these cases, $E^\varphi(\mathbb{R})=L^\varphi(\mathbb{R})$, and the norm and modular convergences are equivalent since the $\Delta_2$-condition is satisfied.

Furthermore, functions such as $\varphi(t)=e^{t^\alpha}-1$, for $\alpha>0$ are $\varphi$-functions do not satisfy the $\Delta_2$-condition, providing cases where $E^\varphi(\mathbb{R}) \neq L^\varphi(\mathbb{R})$, and the modular convergence is weaker than the norm convergence. The Orlicz spaces generated by the latter functions are called \emph{exponential spaces}.

For further details concerning Orlicz spaces see \cite{10, 9, 11, 12, 5, 18}  and \cite{MAJ}.

%%%

\section{Steklov sampling operators and preliminary results} \label{sect3}

In this section, we recall the class of operators we will consider in this paper, recently introduced in \cite{3}. 

These operators aim to reconstruct a given signal $f$ through a series representation built on a kernel function $\chi \colon \mathbb{R} \to \mathbb{R}$ (see \cite{26}, \cite{13}) . The key feature of this approach lies in the use of sample values that are Steklov integrals of order $r$, evaluated at the sampling nodes $\frac{k}{w}$, $k \in \mathbb{Z}$, $w>0$. 
Among the various definitions available in the literature, we adopt a formulation inspired by Sendov and Popov in \cite{1}:
\begin{equation}
f_{r, h}(x):=(-h)^{-r} \int_0^h \cdots \int_0^h \sum_{m=1}^r (-1)^{r-m+1} \binom{r}{m} f \left(x+\frac{m}{r} \left(t_1+ \cdots + t_r \right) \right)\ dt_1 \ldots dt_r, \label{Steklov-type integrals}
\end{equation} 
for any locally integrable function of the form $f \colon \mathbb{R} \to \mathbb{R}$, with $r \in \mathbb{N}^+$, $h>0$.

In what follows, we will use this definition.
\begin{definition} \label{defnucleodiscreto}
A function $\chi \colon \mathbb{R} \to \mathbb{R}$ is called a \textnormal{(}discrete\textnormal{)} kernel if it satisfies the following assumption:
\begin{enumerate}
\item[$(\chi_1)$] $\chi$ is continuous on $\mathbb{R}$;
\item[$(\chi_2)$] the discrete algebraic moment of order $0$: 
\begin{equation}
m_0(\chi, u):= \sum_{k \in \mathbb{Z}} \chi(u-k)=1, \quad u \in \mathbb{R};
\end{equation}
\item[$(\chi_3)$] for some $\alpha>0$, the discrete absolute moment of order $\alpha$:
\begin{equation}
M_\alpha(\chi):= \sup_{u \in \mathbb{R}}\sum_{k \in \mathbb{Z}} |u-k|^\alpha |\chi(u-k)|<+\infty. \label{def nuclei}
\end{equation}
\end{enumerate}
\end{definition}
We generally refer to the assumption $(\chi_2)$ as the partition of the unit property.

The following lemma provides a useful characterization in terms of the Fourier transform of $\chi$ (and its derivatives), which can also be used to verify $(\chi_2)$.

For this reason, we first recall that the usual $L^1$-Fourier transform of $\chi$ is defined as:
$$
\widehat{\chi}(v)\, :=\, \int \chi(x)\, e^{-i x v}\, dx, \quad v \in \R.
$$
We can state what follows.
\begin{lemma} \label{lemmaVMC}
Let $\chi \colon \mathbb{R} \to \mathbb{R}$ be a given function with $\chi \in L^1(\mathbb{R})$, satisfying $(\chi_1)$. In addition, assuming that the function $g(x):= -(ix)^j \chi(x)$, $j \in \mathbb{N}$ \textnormal{(}$x \in \mathbb{R}$ and $i$ denotes the complex unit\textnormal{)} and belongs to $L^1(\mathbb{R})$, the following two assertions are equivalent for every $j \in \N$:
\begin{enumerate}
\item[\textnormal{(a)}]
$$\sum_{k \in \mathbb{Z}}(u-k)^j \ \chi(u-k)=:A_j^\chi \in \mathbb{R}, \quad u \in \mathbb{R};$$
\item[\textnormal{(b)}] $$[\widehat{\chi}]^{(j)}(2k\pi)=\begin{cases}
A_j  & \text{ if } k=0 \\
0 & \text{ if } k \in \mathbb{Z}\setminus \{0\}
\end{cases},$$ where $\widehat{\chi}^{(j)}$ denotes the j-th derivative of $\widehat{\chi}$ (if $j=0$ this condition involves $\widehat{\chi}$),
and  $A_j^\chi=(-i)^jA_j.$
\end{enumerate}
\end{lemma}
A crucial step in the proof of Lemma \ref{lemmaVMC} is played by the celebrated Poisson summation formula, widely used in Fourier Analysis.

Clearly, from Lemma \ref{lemmaVMC} we can deduce an equivalent condition for $(\chi_2)$, that is $\widehat{\chi}(0)=1$ and $\widehat{\chi}(2\pi k)=0$, $k \in \mathbb{Z}\backslash \{0\}$. We can now recall the following.
 
 \begin{definition}
Let $r \in \mathbb{N}^+$ be fixed and let $\chi \in L^1(\mathbb{R})$ be a kernel that satisfies properties $(\chi_1)$, $(\chi_2)$ and $(\chi_3)$ of Definition \textnormal{\ref{defnucleodiscreto}}. We define the Steklov sampling operators of order $r$ $(SSO_r)$, as follows:
\begin{equation}
\begin{split}
&\left(S_w^r f \right)(x):= \sum_{k \in \mathbb{Z}} f_{r, w} \left ( \frac{k}{w} \right) \chi \left(wx-k \right)
=\sum_{k \in \mathbb{Z}} \chi \left( wx-k \right) \times \\ 
& \times \left[ w^r \int_0^\frac{1}{w} \cdots \int_0^\frac{1}{w} \sum_{m=1}^r (-1)^{1-m} \binom{r}{m} f \left(\frac{k}{w}+\frac{m}{r} \left(t_1+ \cdots + t_r \right) \right) \ dt_1 \ldots dt_r \right],
\end{split}
\end{equation}
with $x \in \mathbb{R}$, $w>0$, and the symbol $f_{r, w}$ denotes the Steklov-type integrals defined in \textnormal{(\ref{Steklov-type integrals})} with $h=\frac{1}{w}$, for any locally integrable function $f \colon \mathbb{R} \to \mathbb{R}$ for which the above series are convergent.
 \end{definition}
It is clear that, if we set $r=1$, we obtain $S_w^1=K_w^\chi$, i.e., the Steklov sampling operators coincide with the Kantorovich sampling operators (see, e.g., \cite{2,CNV3}), that are of the form:
$$
(K^\chi_w f)(x)\, :=\, \sum_{k \in \mathbb{Z}} \left[   w \int_{k/w}^{(k+1)/w} f(u)\, du\right]\, \chi(wx-k), \quad x \in \R.
$$
Based on the latter remark, we stress that all the results presented below are significant mainly for the cases $r \geq 2$, where the Steklov sampling operators are very general and different from other known families of sampling type operators.

We now recall the following lemma, which provides useful properties of the discrete kernel $\chi$.

\begin{lemma} \label{lemma}
Under the above assumptions $(\chi_1)$ and $(\chi_3)$ on the kernel $	\chi$, we have:
\begin{enumerate}
\item[\textnormal{(1)}] The discrete absolute moments  
$M_\nu(\chi) < +\infty$, for every $0 \leq \nu \leq \alpha$;
\item[\textnormal{(2)}] for every $\gamma>0$, 
\begin{equation}
\lim_{w \to +\infty} \sum_{|wx-k|>\gamma w} \left | \chi \left(wx-k \right) \right|=0,\label{prop. nuc.}
\end{equation}
uniformly with respect to $x \in \mathbb{R}$;
\item[\textnormal{(3)}] if $\chi \in L^1(\mathbb{R})$, for every $\gamma >0$ and $\varepsilon>0$, there exists a constant $M > 0$ such that: 
$$\int_{|x|>M} w \left | \chi (wx-k)\right | \ dx < \varepsilon,$$ 
for every sufficiently large $w>0$ such that $\frac{k}{w} \in [-\gamma, \ \gamma]$.
\end{enumerate}
\end{lemma}
For a proof of Lemma \ref{lemma}, see, e.g., \cite{2}.

\begin{remark} \label{Note}
\begin{enumerate}
\item[(i)] We note that the previously defined sampling-type series are properly defined for all $r \in \mathbb{N}^+$ and $w > 0$, provided the function $f \colon \mathbb{R} \to \mathbb{R}$ is bounded. Indeed
\begin{equation}
\left| \left(S_w^r f \right)(x) \right| \leq \left( 2^r-1 \right) \| f \|_\infty M_0 (\chi) < +\infty.
\label{absstek} 
\end{equation}
\item[(ii)] Instead of assuming that the function $\chi$ is bounded in a neighborhood of the origin and satisfies $(\chi_3)$, one can directly assume that properties (i) and (ii) of Lemma \textnormal{\ref{lemma}} hold.
\item[(iii)] The moment condition $(\chi_3)$ on the kernels $\chi$ is obviously satisfied if $\chi(x)=\mathcal{O}(x^{-1-\beta-\varepsilon})$ for $x \to \pm \infty$ and some $\varepsilon >0$. Hence, if we consider a kernel $\chi$ with compact support, $(\chi 3)$ is satisfied for every $\alpha>0$.
\end{enumerate}
\end{remark}
%

%%%%

\section{Modular convergence in Orlicz spaces} \label{sect4}

In order to extend the convergence results for the Steklov sampling operators to the general setting of Orlicz spaces generated by convex $\varphi$-functions, we need to recall the following theorem established in \cite{3}. 
\begin{theorem} \label{convStek}
Let $f \colon \mathbb{R} \to \mathbb{R}$ be a bounded function and $r \in \mathbb{N}^+$. Then:
$$\lim_{w \to +\infty} \left(S_w^r f\right)(x)=f(x),$$ at any point of continuity $x \in \mathbb{R}$ of the function $f$. Moreover, in $f \in C(\mathbb{R})$, we have: $$\lim_{w \to +\infty} \| S_w^r f - f \|_\infty=0.$$
\end{theorem}
We can now state the following useful remark.
\begin{remark} \label{remark}
Note that, if $f \in C_C(\mathbb{R})$, with $\text{supp } f \subset [-B, B]$, $B > 0$, and if we choose $\gamma > B + 1$, then for  $w \geq r$ and $k \not \in [- \gamma w, \gamma w]$, we have: \begin{equation}
f_{r, w}\left(\frac{k}{w}\right)=0. \label{eqnull}
\end{equation}
Indeed, for $k<-\gamma w$ and for $t_i \in \left[0, \frac{1}{w}\right]$, where $i=1, \ldots, r$, we have that $\frac{k}{w}<-\gamma$ and
$$\frac{k}{w}+\frac{m}{r}(t_1+ \ldots+t_r) < -\gamma+\frac{m}{w}<-B-1+\frac{m}{w}<-B-1+1=-B.$$
where the last step follows from the fact that $1<m<r \leq w$.
Moreover, if $k>\gamma w$, then $\frac{k}{w}>\gamma$ and since $t_1+ \cdots +t_r \geq 0$, it follows that: 
$$\frac{k}{w}+\frac{m}{r}(t_1+ \cdots + t_r) \geq \gamma>B.$$
\end{remark}
In order to prove a modular convergence theorem for the Steklov sampling operators, we firstly have to show that such operators are well defined in $L^\varphi(\mathbb{R})$ for functions $f \in L^\varphi(\mathbb{R})$, with $\varphi$ convex.

For this reason, from now on, we always consider $\varphi$ as a convex $\varphi$-function and the parameter $r \geq 2$ (since for $r=1$, $S_w^1=K_w$ and hence everything has been already proved in \cite{2}). We begin proving the following proposition.
\begin{proposition} \label{prop1Cap4}
Let $f \in C_c(\mathbb{R})$ be fixed. Then $S_w^r f \in E^\varphi(\mathbb{R})$, $w>0$.
\end{proposition}
\begin{proof}
Let $f \in C_c(\mathbb{R})$ be fixed such that $\text{supp } f \subset [-B, B]$, $B > 0$. As proved in (\ref{eqnull}), if we choose $\gamma > B + 1$, it turns out that for  $w \geq r$ and $k \not \in [- \gamma w, \gamma w]$, we have: 
$$
f_{r, w}\left(\frac{k}{w}\right)=0.$$
Let now $\lambda>0$ be fixed. Using Jensen's inequality and Fubini-Tonelli theorem, we have:
$$I^\varphi \left[\lambda S_w^r f \right]=$$
$$= \int_\mathbb{R} \varphi \left(\lambda \left| \left(S_w^r f \right)(x)\right| \right) \ dx= \int_\mathbb{R} \varphi \left(\lambda \left| \sum_{k \in \mathbb{Z}} f_{r, w} \left(\frac{k}{w}\right) \chi(wx-k) \right| \right) \ dx $$
$$= \int_\mathbb{R} \varphi \left(\lambda \left| \sum_{|k| \leq \gamma w} f_{r, w} \left(\frac{k}{w}\right) \chi(wx-k) \right| \right) \ dx$$
$$ \leq  \int_\mathbb{R} \varphi \left(\lambda  \sum_{|k| \leq \gamma w} \left| f_{r, w}  \left(\frac{k}{w}\right)\right| \left| \chi(wx-k) \right| \right) \ dx $$
$$\leq \sum_{|k| \leq \gamma w} \frac{\varphi \left(\lambda M_0(\chi) \left|f_{r, w} \left(\frac{k}{w}\right) \right| \right)}{M_0(\chi)} \int_{\mathbb{R}} \left|\chi(wx-k)\right| \ dx.$$
Using the change of variables $u=wx-k$, we obtain:
$$I^\varphi \left[\lambda S_w^r f \right] \leq \sum_{|k| \leq \gamma w} \frac{\varphi \left(\lambda M_0(\chi) \left|f_{r, w} \left(\frac{k}{w}\right) \right| \right)}{M_0(\chi)} \cdot \frac{1}{w}\int_{\mathbb{R}} \left|\chi(u)\right| \ du $$
$$= \sum_{|k| \leq \gamma w} \frac{\varphi \left(\lambda M_0(\chi) \left|f_{r, w} \left(\frac{k}{w}\right) \right| \right)}{w M_0(\chi)} \|\chi\|_1.$$
Now, using Jensen's inequality again, we can analyze the following expression: 
\setlength{\jot}{15pt}
$$\left| f_{r, w} \left( \frac{k}{w}\right) \right|=$$
$$=\left| w^r \int_0^{1/w} \cdots \int_0^{1/w} \sum_{m=1}^r (-1)^{1-m} \binom{r}{m} f\left(\frac{k}{w}+\frac{m}{r}(t_1 + \cdots +t_r) \right) \ dt_1 \ldots dt_r \right|$$
$$\leq w^r \int_0^{1/w} \cdots \int_0^{1/w} \left|   \sum_{m=1}^r (-1)^{1-m} \binom{r}{m} f\left(\frac{k}{w}+\frac{m}{r}(t_1 + \cdots +t_r) \right) \right| \ dt_1 \ldots dt_r$$
$$\leq w^r \int_0^{1/w} \cdots \int_0^{1/w} \sum_{m=1}^r \binom{r}{m} \left| f \left(\frac{k}{w}+\frac{m}{r}(t_1 + \cdots +t_r) \right) \right| \ dt_1 \ldots dt_r $$
\begin{equation}
\leq (2^r-1) w^r \cdot \frac{1}{w^r} \|f\|_\infty=(2^r-1) \|f\|_\infty. \label{disprop1Orlicz}
\end{equation}

In addition, we can also compute the sum: 
\begin{equation}
\sum_{|k|<\gamma w} \frac{1}{w} \leq 2 \sum_{k=0}^{\lfloor \gamma w \rfloor} \frac{1}{w}=2\left[\lfloor \gamma w \rfloor+1 \right]\cdot \frac{1}{w} \leq 2\gamma +1, \label{eqsumw}
\end{equation} 
where $\lfloor \boldsymbol{\cdot} \rfloor$ denotes the integer part of the given real number.

Finally, using (\ref{eqsumw}), we have:
$$I^\varphi \left[\lambda S_w^r f \right] \leq \sum_{|k| \leq \gamma w} \frac{\varphi \left(\lambda M_0(\chi) (2^r-1) \|f\|_\infty \right)}{w M_0(\chi)} \|\chi\|_1$$ 
$$\leq \frac{\varphi \left(\lambda M_0(\chi) (2^r-1) \|f\|_\infty \right)}{ M_0(\chi)} (2 \gamma+1)\|\chi\|_1 < +\infty.$$
This completes the proof.
\end{proof}

We can now prove the following norm convergence theorem for the family of Steklov sampling operators in Orlicz spaces, when $f \in C_C(\mathbb{R})$.

\begin{theorem} \label{thnormphi}
For every $f \in C_c(\mathbb{R})$, we have $$\lim_{w \to +\infty} \|S_w^r f-f\|_\varphi=0.$$ 
\end{theorem}

\begin{proof}
First, we can note that, by Proposition \ref{prop1Cap4} and since $C_C(\mathbb{R}) \subset E^\varphi(\mathbb{R})$, it turns out that $f$, $S_w^r f \in E^\varphi(\mathbb{R})$, $w>0$. Moreover, to prove the thesis it is sufficient to show that
$$\lim_{w \to +\infty} I^\varphi\left[\lambda(S_w^r f -f)\right]=\lim_{w \to +\infty} \int_\mathbb{R} \varphi \left(\lambda |(S_w^r f)(x)-f(x)|\right) \ dx=0,$$ 
for every $\lambda>0$, that is equivalent to show that the family $\varphi\left(\lambda \left|S_w^r f-f\right|\right)_{w>0}$ converges to zero in $L^1(\mathbb{R})$, for every $\lambda>0$. To show the latter fact, we will use the well known Vitali convergence theorem (VCT) in $L^1(\mathbb{R})$.

Let $\lambda>0$ be fixed. Using Theorem \ref{convStek} and by the continuity of $\varphi$, we have
$$\lim_{w \to +\infty}\varphi\left(\lambda \|S_w^r f-f\|_\infty \right)=0.$$
This implies that the first condition of the VCT is satisfied, i.e., the considered sequence converges in measure to zero.
Now, in order to prove the second one, let $\varepsilon>0$ be fixed. If $\text{supp } f \subset [-B, B]$, $B > 0$, and if
we choose $\gamma > B + 1$, then for  $w \geq r$ and $k \not \in [- \gamma w, \gamma w]$, we have: $$
f_{r, w}\left(\frac{k}{w}\right)=0,$$
as showed in Remark \ref{remark}.

Now, by Lemma \ref{lemma} (3), there exists a constant $M>0$ (we can assume $M>B$), such that
$$\int_{|x|>M} w |\chi(wx-k)| \ dx <\varepsilon,$$
for $w>0$ sufficiently large and $k \in [-\gamma w, \gamma w]$.

Then, by Jensen's inequality and Fubini-Tonelli theorem, we have
$$\int_{|x|>M} \varphi \left( \lambda \left| \left(S_w^r f\right)(x)\right|\right) \ dx=\int_{|x|>M} \varphi \left( \lambda \left| \sum_{k \in \mathbb{Z}} f_{r, w}\left(\frac{k}{w}\right) \chi(wx-k)\right|\right) \ dx$$
$$=\int_{|x|>M} \varphi \left( \lambda  \sum_{k \in \mathbb{Z}} \left| f_{r, w}\left(\frac{k}{w}\right)\right| \left| \chi(wx-k)\right|\right) \ dx $$ 
$$\leq \sum_{k \in [-\gamma w, \gamma w]} \frac{1}{w M_0(\chi)} \varphi\left(\lambda M_0(\chi) \left|f_{r, w} \left(\frac{k}{w}\right)\right| \right) \int_{|x|>M} w |\chi(wx-k)| \ dx$$ 
$$\leq \frac{\varepsilon}{w M_0(\chi)} \sum_{k \in [-\gamma w, \gamma w]}\varphi \left(\lambda M_0(\chi) \left|f_{r, w} \left(\frac{k}{w}\right) \right| \right)$$ 
$$\leq \frac{\varepsilon}{w M_0(\chi)} \sum_{k \in [-\gamma w, \gamma w]}\varphi \left(\lambda M_0(\chi) (2^r-1) \|f\|_\infty \right),$$
where in the above inequality we used (\ref{disprop1Orlicz}) again.

Now, noting that the number of terms in the above series in the foregoing inequality
does not exceed $2\left(\lfloor\gamma  w\rfloor+1\right)$ (see  (\ref{eqsumw}) again), so we have
$$\int_{|x|>M} \varphi \left(\lambda \left| \left(S_w^r f\right)(x)\right|\right) \ dx \leq \frac{2\varepsilon}{M_0(\chi)} \left(\gamma +1\right) \varphi\left(\lambda M_0(\chi)(2^r-1) \|f\|
_\infty\right)=:\varepsilon \cdot C.$$
Therefore, we proved that for $\varepsilon>0$, there exists a set $E_\varepsilon=[-M, M]$ such that, for every measurable set $F$, with $F \cap E_\varepsilon=\emptyset$, we have
$$\int_F \varphi\left(\lambda\left|\left(S_w^r f\right)(x)-f(x)\right|\right) \ dx= \int_F \varphi\left(\lambda \left| \left(S_w^r f \right)(x)\right| \right) \ dx $$
$$\leq \int_{|x|>M } \varphi\left(\lambda\left| \left(S_w^r f \right)(x) \right| \right) \ dx < \varepsilon \cdot C,$$
namely the second condition of the VCT holds.

Finally, we prove the last condition of the VCT, i.e., that all the integrals of $\varphi\left(\lambda |(S_w^r f)(x)-f(x)|\right)$ are equi-absolutely continuous. For $\|f\|_\infty >0$, let $\varepsilon>0$ and let $B_\varepsilon \subset \mathbb{R}$  a measurable set with measure less than 
$$\delta:=\frac{\varepsilon}{\varphi \left(2 \lambda (2^r-1)M_0(\chi)\|f\|_\infty\right)}.$$ 
Hence, using (\ref{absstek}), the convexity of $\varphi$, $M_0(\chi) \geq 1$ and $2^r-1 \geq 1$, we have
$$ \int_{B_\varepsilon} \varphi \left(\lambda\left|(S_w^r f)(x)-f(x) \right| \right) \ dx \leq \int_{B_\varepsilon} \varphi \left( \lambda \left| \left(S_w^r f\right)(x) \right| + \left|f(x) \right| \right) \ dx$$ 
$$\leq \frac{1}{2} \int_{B_\varepsilon} \varphi \left(2 \lambda \left| \left(S_w^r f\right)(x) \right| \right) \ dx + \frac{1}{2}\int_{B_\varepsilon} \varphi\left(2\lambda \left|f(x) \right| \right) \ dx $$ 
$$\leq \frac{1}{2} \int_{B_\varepsilon} \varphi \left(2 \lambda (2^r-1) \|f\|_\infty M_0(\chi) \right) \ dx + \frac{1}{2} \int_{B_\varepsilon} \varphi \left(2\lambda\|f\|_\infty \right) \ dx$$
 $$\leq \int_{B_\varepsilon} \varphi\left(2\lambda(2^r-1)M_0(\chi)\|f\|_\infty \right) \ dx < \varphi\left(2\lambda(2^r-1)M_0(\chi) \|f\|_\infty \right) \delta < \varepsilon.$$
Now, since $\lambda>0$ is arbitrary, we obtain the assertion. 
 \end{proof}
In order to obtain a modular convergence result in $L^\varphi(\mathbb{R})$, we need a modular continuity property for the operators $S_w^r$. We can prove the following.

\begin{theorem}\label{modconv}
 For every $f \in L^\varphi(\mathbb{R})$, there holds:
$$
I^\varphi[\lambda S_w^r f] \leq \frac{\| \chi\|_1}{M_0(\chi)}r I^\varphi\left[\lambda (2^r-1) M_0(\chi) f \right], \quad w>0
$$
for some $\lambda>0$.

In particular, $S_w^r f$ is well defined and $S_w^r f \in L^\varphi(\mathbb{R})$ whenever $f \in L^\varphi(\mathbb{R})$.
\end{theorem}
\begin{proof}
Since $f \in L^\varphi(\mathbb{R})$, there exists $\overline{\lambda}>0$ such that $I^\varphi\left[ \overline{\lambda}f\right]<+\infty$. Let $\lambda>0$ such that $\lambda (2^r-1) M_0(\chi) \leq \overline{\lambda}$, then applying  Jensen's inequality twice, and the change of variable $u=wx-k$, we obtain
$$
I^\varphi \left[\lambda S_w^r f \right]= \int_\mathbb{R} \varphi \left(\lambda \left| \left(S_w^r f \right)(x)\right| \right) \ dx= \int_\mathbb{R} \varphi \left(\lambda \left| \sum_{k \in \mathbb{Z}} f_{r, w} \left(\frac{k}{w}\right) \chi(wx-k) \right| \right) \ dx$$
$$\leq \int_\mathbb{R} \varphi \left(\lambda \sum_{k \in \mathbb{Z}}  \left| f_{r, w} \left(\frac{k}{w}\right)\right| \left| \chi(wx-k) \right| \right) \ dx  $$
$$\leq\frac{1}{M_0(\chi)} \sum_{k \in \mathbb{Z}} \varphi \left(\lambda M_0(\chi) \left| f_{r, w}  \left(\frac{k}{w}\right)\right|\right) \int_\mathbb{R}  \left| \chi(wx-k) \right| \ dx $$
$$= \frac{1}{M_0(\chi)} \sum_{k \in \mathbb{Z}} \varphi \left(\lambda M_0(\chi) \left|f_{r, w} \left(\frac{k}{w}\right) \right| \right) \cdot \frac{1}{w}\int_{\mathbb{R}} \left|\chi(u)\right| \ du$$
$$= \frac{\|\chi\|_1}{w M_0(\chi)} \sum_{k \in \mathbb{Z}} \varphi \left(\lambda M_0(\chi) \left|f_{r, w} \left(\frac{k}{w}\right) \right| \right).$$
Using Jensen's inequality again and the change of variable $y=\frac{m}{r}t_1+\frac{k}{w}$, we have:
$$\varphi \left( \lambda M_0(\chi) \left| f_{r, w} \left(\frac{k}{w}\right) \right| \right) \ \leq $$
$$\leq \varphi \left( \lambda M_0(\chi) w^r \sum_{m=1}^r \binom{r}{m} \int_0^{1/w} \cdots \int_0^{1/w}  \left| f \left(\frac{k}{w}+\frac{m}{r}(t_1 + \cdots +t_r) \right) \right| \ dt_1 \ldots dt_r \right)$$
$$\leq \frac{1}{2^r-1} \cdot w^r \sum_{m=1}^r \binom{r}{m} \int_0^{1/w} \cdots \int_0^{1/w} \varphi \biggl( \lambda (2^r-1) M_0(\chi)$$
$$\left| f \left(\frac{k}{w}+\frac{m}{r}(t_1 + \cdots +t_r) \right) \right| \biggr) \ dt_1 \ldots dt_r $$
$$ = \frac{1}{2^r-1} \cdot w^r \sum_{m=1}^r \binom{r}{m}\int_0^{1/w} \cdots \int_0^{1/w}  dt_2 \ldots dt_r$$ 
$$\frac{r}{m} \int_{\frac{k}{w}}^{\frac{k}{w}+\frac{m}{rw}} \varphi \left( \lambda (2^r-1) M_0(\chi) \left| f \left(y+\frac{m}{r}(t_2 + \cdots +t_r) \right) \right| \right) \ dy$$
$$ \leq \frac{1}{2^r-1} \cdot w^r r \sum_{m=1}^r \binom{r}{m} \int_0^{1/w} \cdots \int_0^{1/w}  dt_2 \ldots dt_r$$
$$ \int_{\frac{k}{w}}^{\frac{k+1}{w}} \varphi \left(\lambda (2^r-1) M_0(\chi) \left| f \left(y+\frac{m}{r}(t_2 + \cdots +t_r) \right) \right| \right) \ dy.$$
Finally, returning to the previous estimate, we have:
$$I^\varphi \left[\lambda S_w^r f \right] \leq$$ 
$$\leq \frac{\|\chi\|_1}{w M_0(\chi)} \sum_{k \in \mathbb{Z}} \biggl[\frac{1}{2^r-1} \cdot w^r r \sum_{m=1}^r \binom{r}{m} \int_0^{1/w} \cdots \int_0^{1/w}  dt_2 \ldots dt_r$$
$$ \int_{\frac{k}{w}}^{\frac{k+1}{w}} \varphi \left(\lambda (2^r-1) M_0(\chi) \left| f \left(y+\frac{m}{r}(t_2 + \cdots +t_r) \right) \right| \right) \ dy \biggr]$$
$$\leq \frac{1}{2^r-1} \cdot\frac{\|\chi\|_1}{w M_0(\chi)} w^r r \sum_{m=1}^r \binom{r}{m} \int_0^{1/w} \cdots \int_0^{1/w}  dt_2 \ldots dt_r$$
$$\sum_{k \in \mathbb{Z}} \int_{\frac{k}{w}}^{\frac{k+1}{w}} \varphi \left(\lambda (2^r-1) M_0(\chi) \left| f \left(y+\frac{m}{r}(t_2 + \cdots +t_r) \right) \right| \right) \ dy$$
$$\leq \frac{1}{2^r-1} \cdot \frac{\|\chi\|_1}{M_0(\chi)} w^{r-1} r \sum_{m=1}^r \binom{r}{m} \int_0^{1/w} \cdots \int_0^{1/w}  dt_2 \ldots dt_r$$
$$I^\varphi \left[\lambda (2^r-1) M_0(\chi) f\left(\boldsymbol{\cdot}+\frac{m}{r}(t_2+ \cdots t_r)\right)\right]$$
$$= \frac{\|\chi\|_1}{M_0(\chi)} w^{r-1} r \cdot \frac{1}{w^{r-1}} I^\varphi \left[\lambda(2^r-1)M_0(\chi) f(\boldsymbol{\cdot})\right]$$
$$= \frac{\|\chi\|_1}{M_0(\chi)} r I^\varphi \left[\lambda (2^r-1) M_0(\chi) f \right] < +\infty,$$
since $I^\varphi \left[\lambda(2^r-1) M_0(\chi)f\right] \leq I^\varphi \left[\overline{\lambda}f\right]<+\infty$ and we know that any $$I^\varphi \left[\lambda (2^r-1) M_0(\chi) f\left(\boldsymbol{\cdot}+\frac{m}{r}(t_2+ \cdots t_r)\right)\right]= I^\varphi \left[\lambda (2^r-1) M_0(\chi) f\left(\boldsymbol{\cdot}\right)\right].$$
This completes the proof.
\end{proof}

Theorem \ref{modconv} says that the family of operators $S_w^r \colon L^\varphi(\mathbb{R}) \to L^\varphi(\mathbb{R})$ is modularly continuous, for every $w >0$.

We can finally prove the main theorem of this section.
\begin{theorem}
For every $f \in L^\varphi(\mathbb{R})$, there exists $\lambda>0$ such that
$$\lim_{w \to +\infty} I^\varphi \left[\lambda \left(S_w^r f -f\right)\right]=0.$$
\end{theorem}
\begin{proof}
Let $f \in L^\varphi(\mathbb{R})$. By the modularly density property of $C_C(\mathbb{R})$ in $L^\varphi(\mathbb{R})$, there exists $\overline{\lambda}>0$ such that, for every $\varepsilon>0$, there exists a function $g \in C_c(\mathbb{R})$ such that $I^\varphi \left[ \overline{\lambda}(f-g) \right] < \varepsilon$. 

Now, let $\lambda>0$ be such that $3\lambda\left(1+(2^r-1) M_0(\chi)\right) \leq \overline{\lambda}$.
By the properties of $\varphi$ and Theorem \ref{modconv}, we have
$$I^\varphi \left[\lambda\left(S_w^r f-f \right) \right] \leq $$
$$ \leq I^\varphi \left[3\lambda\left(S_w^r f- S_w^r g \right) \right]+I^\varphi \left[3\lambda\left(S_w^r g- g \right) \right]+I^\varphi \left[3\lambda\left( f- g \right) \right]$$
$$\leq \frac{\|\chi\|_1 r}{M_0(\chi)}I^\varphi\left[\overline{\lambda}\left(f-g\right)\right]+ I^\varphi \left[3\lambda\left(S_w^r g- g \right) \right]+I^\varphi \left[ \overline{\lambda} \left( f- g \right) \right]$$
$$\leq \left(\frac{\|\chi\|_1 r}{M_0(\chi)}+1\right) \varepsilon + I^\varphi \left[3\lambda\left(S_w^r g- g \right) \right].$$
The assertion now follows from Theorem \ref{thnormphi}. This completes the proof.
\end{proof}
%

%%%%

\section{Particular cases and examples} \label{Particular cases}

In this section, we provide some significant cases in which the previous results can be applied. First, we can note that, if we assume $\varphi(u)=u^p$, for $u \in \mathbb{R}_0^+$, we immediately obtain the convergence results established in \cite{3} in the $L^p$-setting.

As another particular case, we can consider the Orlicz spaces generated by the $\varphi$-function $\varphi_{\alpha, \beta}(u)=u^\alpha \log^\beta(e+u)$, $u \geq 0$ for $\alpha \geq 1$ and $\beta>0$.
The associated Orlicz space consists by the set of all functions $f \in M(\mathbb{R})$ that satisfy the following condition:
$$I^{\varphi_{\alpha, \beta}}[\lambda f]= \int_\mathbb{R} \left(\lambda |f(x)|\right)^\alpha \log^\beta\left(e+\lambda |f(x)|\right) \ dx < +\infty,$$
for some $\lambda>0$, and it is denoted by $L^\alpha \log^\beta L(\mathbb{R})$.
Note that the function $\varphi_{\alpha, \beta}$ satisfies the $\Delta_2$-condition, which means that, as in the case of $L^p$ spaces, we have $L^\alpha \log^\beta L(\mathbb{R})=L^{\varphi_{\alpha, \beta}}(\mathbb{R})=E^{\varphi_{\alpha, \beta}}(\mathbb{R})$, and the norm convergence is equivalent to the modular convergence.

As a result of general theory, we can state the following corollary, e.g., for $\alpha=\beta=1$. 

\begin{corollary} \label{corollario5.1}
For every $f \in L\log L(\mathbb{R})$ and for every $\lambda>0$, there hods:
$$\lim_{w \to +\infty} \int_\mathbb{R} \left|(S_w^r f)(x)-f(x)\right| \log\left(e+\lambda \left|(S_w^r f)(x)-f(x)\right|\right) \ dx=0,$$
or equivalently, 
$$\lim_{w \to +\infty} \|S_w^r f-f\|_{L\log L}=0,$$
where $\| \boldsymbol{\cdot}\|_{L \log L}$ is the Luxemburg norm associated to $I^{\varphi_{1, 1}}$.

Moreover, we have
$$\int_\mathbb{R} \left| S_w^r f(x)\right| \log \left(e+\lambda \left|S_w^r f(x) \right| \right) \ dx \leq$$
$$\leq \|\chi\|_1 r (2^r-1) \int_\mathbb{R} |f(x)| \log\left(e+\lambda (2^r-1) r M_0(\chi) |f(x)|\right) \ dx,$$
for every $\lambda>0$, i.e., $S_w^r \colon L \log L(\mathbb{R}) \to L \log L(\mathbb{R})$.
\end{corollary}

Finally, we consider the case of exponential spaces generated by the $\varphi$-function $\varphi_\alpha(u)=e^{u^\alpha}-1$, $u \geq 0$ for some $\alpha>0$. The Orlicz space $L^{\varphi_\alpha}(\mathbb{R})$ consists of those functions $f \in M(\mathbb{R})$ for which 
$$I^{\varphi_\alpha}[\lambda f]=\int_\mathbb{R} \left(e^{(\lambda |f(x)|)^\alpha}-1\right) \ dx < +\infty,$$
for some $\lambda>0$. Since $\varphi_\alpha$ does not satisfy the $\Delta_2$-condition, the Orlicz space $L^{\varphi_\alpha}(\mathbb{R})$ does not coincide with the space of its finite elements $E^{\varphi_\alpha}(\mathbb{R})$. Furthermore, modular convergence does not necessarily lead to norm convergence.

\begin{corollary} \label{corollario5.2}
For $f \in L^{\varphi_\alpha}(\mathbb{R})$, there exists $\lambda>0$ such that
$$\lim_{w \to +\infty} \int_\mathbb{R} \left(e^{\left(\lambda \left| \left(S_w^r f\right)(x)-f(x) \right| \right)^\alpha}-1\right) \ dx=0.$$

Moreover, we also have
$$\int_\mathbb{R} \left(e^{\left(\lambda \left| \left(S_w^r f\right)(x)-f(x) \right| \right)^\alpha}-1\right) \ dx \leq
\frac{\|\chi\|_1 r}{M_0(\chi)}
  \int_\mathbb{R} \left(e^{(\lambda (2^r-1) M_0(\chi)|f(x)|)^\alpha}-1\right) \ dx.
$$
\end{corollary}
%

%%%%

\section{Applications and further examples}  \label{sec6}

As a first particular case, let us consider the Steklov sampling series based on the Fejér's kernel $F$. The function $F$ is defined as follows 
$$
F(x) := \frac{1}{2} \text{sinc}^2 \left(\frac{x}{2}\right), \quad x \in \mathbb{R},
$$
where the sinc-function and its Fourier transform (in the $L^2$-sense) are given by
$$
\text{sinc} (x) := 
\begin{cases}
\frac{\sin \pi x}{\pi x}, & x \in \mathbb{R} \setminus \{0\}, \\
1, & x = 0,
\end{cases}
\qquad
\widehat{\text{sinc}} (v) =
\begin{cases}
1, \quad |v| \leq \pi, \\
0, \quad |v| > \pi.
\end{cases}
$$

It is well known that $F \in L^1(\mathbb{R})$ and is bounded on $\R$. Furthermore, it satisfies condition $(\chi_3)$ with any $0<\beta<1$ in view of Remark \ref{Note}(iii) (while $M_1(F)=+\infty$, see \cite{CACO2025}). Since its Fourier transform is given by (see, e.g., \cite{26})
$$\widehat{F}(v)=
\begin{cases}
    1-\left|\frac{v}{\pi}\right|, \quad v \leq \pi\\
    0, \quad \quad \quad \quad v>\pi,
\end{cases}$$
it follows from Remark \ref{Note}(iv) that $F$ satisfies the partition of the unit property $(\chi_2)$.

The Steklov sampling operators, for $f \in L^p(\mathbb{R})$, $1 \leq p < +\infty$, are now given by
$$
\left(S_w^{r,F} f \right)(x):= \sum_{k \in \mathbb{Z}} \frac{1}{2} \ \text{sinc}^2\left(\frac{wx-k}{2}\right) \times $$ 
$$ \times \left[ w^r \int_0^\frac{1}{w} \cdots \int_0^\frac{1}{w} \sum_{m=1}^r (-1)^{1-m} \binom{r}{m} f \left(\frac{k}{w}+\frac{m}{r} \left(t_1+ \cdots + t_r \right) \right) \ dt_1 \ldots dt_r \right].
$$
From $L^p$-convergence theorem for functions belonging to $L^p$ (see Theorem 4.4 in \cite{3}), Corollary \ref{corollario5.1} and Corollary \ref{corollario5.2}, we immediately have the following.

\begin{corollary} \label{corol}
\begin{enumerate}[label=\textnormal{\arabic*.}]
 \item For every $f \in L^p(\mathbb{R})$, $1 \leq p < \infty$, we have
$$\lim_{w \to +\infty} \left\| S_w^{r,F} f - f \right\|_p = 0.$$
\item For $f \in L \log L$ and for every $\lambda > 0$, we have
$$\lim_{w \to +\infty} \int_{\mathbb{R}} \left| S_w^{r,F} f(x) - f(x) \right| \log\left( e + \lambda \left| S_w^{r,F} f(x) - f(x) \right| \right) dx = 0$$
or, equivalently,
$$\lim_{w \to +\infty} \left\| S_w^{r,F} f - f \right\|_{L \log L} = 0.
$$
\item For $f \in L^{\varphi_\alpha}(\mathbb{R})$, there exists $\lambda>0$ such that $$\lim_{w \to +\infty} \int_{\mathbb{R}} \left( e^{\left( \lambda \left| S_w^{r,F} f(x) - f(x) \right| \right)^\alpha} - 1 \right) dx = 0.
$$
\end{enumerate}
\end{corollary}

We note that the Fejér's kernel represents an example of (band-limited) kernel with unbounded support, characterized by a decay of the type ${\cal O}(1/x^2)$, as $x \to \pm \infty$. 

One can also consider kernels $\chi(x)$ which have a faster decay, such as the Jackson type kernels
$$
J_n(x) = c_n  \text{sinc}^{2n} \left( \frac{x}{2n\pi\alpha} \right), \qquad \text{with } \
 c_n := \left[ \int_{-\infty}^{+\infty} \text{sinc}^{2n} \left( \frac{u}{2n\pi\alpha} \right) du \right]^{-1},
$$
$x \in \mathbb{R}$, $n \in \mathbb{N}$ and $\alpha \geq 1$. The $J_n$ are again examples of kernel that are band-limited in $[-1/\alpha, 1/\alpha]$, i.e., their Fourier transform vanishes outside this interval, condition $(\chi_2)$ is satisfied, and $(\chi_3)$ holds in view of 
Remark \ref{Note} (iii), (iv), since $J_n$ decays as ${\cal O}(x^{-2n})$, as $x \to \pm \infty$. Thus, Corollary \ref{corol} holds again if we replace the Fejér's kernel $F$ by the Jackson type kernels $J_n$.

Sometimes, in real-world or in numerical applications it can be useful to deal with kernels with compact support (duration limited).
The latter fact can help us to reduce certain numerical errors, such as the typical truncation error characterizing algorithms theoretically based on series.

Among the most useful examples of kernels with compact support we can mention the central B-splines of order
$n \in \mathbb{N}$, defined by
$$
M_n(x) := \frac{1}{(n - 1)!} \sum_{j=0}^{n} (-1)^j \binom{n}{j} \left( \frac{n}{2} + x - j \right)_+^{n - 1},
$$
where $(x)^r_+:=\max \{x^r, 0\}$,
with Fourier transform given by 
$$
\widehat{M}_n(v) = \text{sinc}^n \left( \frac{v}{2\pi} \right),  \quad v \in \mathbb{R}.
$$

\noindent The $M_n$ are bounded in $\mathbb{R}$ for all $n \in \mathbb{N}$ with compact support $[-n/2, n/2]$.
This implies that $M_n \in L^1(\mathbb{R})$ and that the condition $(\chi_3)$ is satisfied
for all $\beta > 0$. Finally, condition $(\chi_2)$ holds again in view of
Remark \ref{Note} (iv). So, we have
$$
\left(S_w^{r,M_n} f \right)(x):= \sum_{k \in \mathbb{Z}} M_n(wx-k) \times $$ 
$$ \times \left[ w^r \int_0^\frac{1}{w} \cdots \int_0^\frac{1}{w} \sum_{m=1}^r (-1)^{1-m} \binom{r}{m} f \left(\frac{k}{w}+\frac{m}{r} \left(t_1+ \cdots + t_r \right) \right) \ dt_1 \ldots dt_r \right].
$$
Hence, also in this case, an analogous result to Corollary \ref{corol} holds.

%
%%%%%%%

\section{Future works and open problems} \label{sec7}

In the paper \cite{3} the author introduced Steklov sampling operators showing that they own better approximation properties (in terms of the order of approximation) than the sampling type series known in the literature. Here, we established more general convergence results in the very wide setting of Orlicz spaces, without considering the problem of establishing estimates for the approximation errors. Thus, in a future work we aim to face the latter open problem; to do this we will consider a suitable definition of the modulus of smoothness in Orlicz spaces. It is known that (see \cite{COPI2}) in Orlicz spaces we have two different notions of modulus of smoothness: the first one defined through the Luxemburg norm (which implies the norm convergence) and a second one based on the modular $I^\varphi$ from which one can deduce the degree of modular convergence. At the same time, from the two moduli, two definitions of Lipschitz classes arise (the strong and weak ones) that are not equivalent in general. From these considerations, we aim to estimate the asymptotic qualitative error of approximation.
On this analysis (we expect that) a crucial role will be played by the assumption of a suitable Strang-Fix type condition related to polynomials preservation.

\section*{Acknowledgments}

{\small The first author is member of the Gruppo Nazionale per l'Analisi Matematica, la Probabilit\`a e le loro Applicazioni (GNAMPA) of the Istituto Nazionale di Alta Matematica (INdAM), of the network RITA (Research ITalian network on Approximation), and of the UMI (Unione Matematica Italiana) group T.A.A. (Teoria dell'Approssimazione e Applicazioni). 
}

\section*{Funding}

{\small The first author has been supported within the project PRIN 2022: ``AI- and DIP-Enhanced DAta Augmentation for Remote Sensing of Soil Moisture and Forest Biomass (AIDA)" funded by the European Union under the Italian National Recovery and Resilience Plan (NRRP) of NextGenerationEU, under the Italian Ministry of University and Research (MUR) (Project Code: 20229FX3B9, CUP: J53D23000660001). 

}
\section*{Conflict of interest/Competing interests}

{\small The author declares that he has no conflict of interest and competing interest.}

\section*{Availability of data and material and Code availability}

{ \small Not applicable.}

%

%%%%%%
\end{document}